\documentclass{IEEEtran}
\usepackage{cite}
\usepackage{amsmath,amssymb,amsfonts}
\usepackage{textcomp}
\usepackage{xcolor}
\usepackage{graphicx}
\usepackage{subfig}
\usepackage{epstopdf}
\usepackage{comment}
\usepackage{amsthm}
\usepackage{algorithmic}
\usepackage{textcomp}
 
\usepackage{etoolbox}
\newtheoremstyle{mystl}
  {0} 
  {0} 
  {\itshape}
  {}
  {\bfseries}  
  {.} 
  {5pt plus 1pt minus 1pt} 
 {\thmname{#1} \thmnumber{#2}\ifblank{#3}{}{ (\thmnote{#3})}} 

\theoremstyle{mystl}
\newtheorem{lem}{Lemma}
\newtheorem{thme}{Theorem}

\newtheorem{prop}{Proposition}
\newtheorem{ass}{Assumption}

\newtheorem{rmk}{Remark}

\renewcommand{\qed}{\hfill\blacksquare}
\def\BibTeX{{\rm B\kern-.05em{\sc i\kern-.025em b}\kern-.08em
    T\kern-.1667em\lower.7ex\hbox{E}\kern-.125emX}}
\begin{document}
\title{Observer-Based Periodic Event-Triggered and Self-Triggered Boundary Control of a Class of Parabolic PDEs}
\author{Bhathiya Rathnayake, \IEEEmembership{Student Member, IEEE} and Mamadou Diagne, \IEEEmembership{Member, IEEE}
\thanks{B. Rathnayake is with the Department of Electrical and Computer Engineering, University of California San Diego, 9500 Gilman Dr, La Jolla, CA 92093. Email: brm222@ucsd.edu}
\thanks{M. Diagne is with the Department of Mechanical and Aerospace Engineering, University of California San Diego, 9500 Gilman Dr, La Jolla, CA 92093. Email: mdiagne@ucsd.edu}
\thanks{This work is supported by the NSF CAREER Award CMMI-2302030 and the NSF grant CMMI-2222250}
}

\maketitle

\begin{abstract}
This paper introduces the first observer-based periodic event-triggered control (PETC) and self-triggered control (STC) for boundary control of a class of parabolic PDEs using PDE backstepping control. We introduce techniques to convert a certain class of continuous-time event-triggered control into PETC and STC, eliminating the need for continuous evaluation of the triggering function. For the PETC, the triggering function requires only periodic evaluations to detect events, while the STC proactively computes the time of the next event right at the current event time using the system model and the continuously available measurements. For both strategies, the control input is updated exclusively at events and is maintained using a zero-order hold between events. We demonstrate that the closed-loop system is Zeno-free. We offer criteria for selecting an appropriate sampling period for the PETC and for determining the time until the next event under the STC. We prove the system's global exponential convergence to zero in the spatial $L^2$ norm for both anti-collocated and collocated sensing and actuation under the PETC. For the STC, local exponential convergence to zero in the spatial $L^2$ norm for collocated sensing and actuation is proven. Simulations are provided to illustrate the theoretical claims.
\end{abstract}


\section{Introduction}
Event-triggered control (ETC) updates the control input based on events generated by a suitable triggering mechanism instead of at fixed intervals. This approach incorporates feedback into the control update tasks, allowing the control input to be updated aperiodically and only when necessary, based on the system's states. In ETC, a primary challenge is avoiding Zeno behavior—infinite updates in a finite time—which is usually achieved by the careful design of the triggering mechanism so that it is endowed with a positive lower bound for the time between events, known as the \textit{minimal dwell-time (MDT)}. Research over the last decade has expanded from ODEs \cite{heemels2012introduction} to those on PDEs, generating considerable advances (e.g. \cite{espitia2020observer,espitia2019event,katz2020boundary,diagne2021event,rathnayake2021observer,rathnayake2022sampled,wang2022event,rathnayake2022observer}). Particularly, \cite{rathnayake2021observer} and \cite{rathnayake2022sampled} are the most relevant to the present study.

A significant limitation of ETC for both ODEs and PDEs is the necessity for continuous-time evaluation of the triggering function to detect events, which is not ideal for digital implementations. These strategies are referred to as continuous-time ETC (CETC). To address this limitation, two alternative approaches have been developed: \textit{periodic event-triggered control (PETC)} which checks the event-triggering function periodically and decides on control updates \cite{heemels2012periodic}, and \textit{self-triggered control (STC)} which proactively calculates the next event time at the current event time, using system states and the knowledge of the system's dynamics \cite{heemels2012introduction}.
Recent works on both PETC \cite{heemels2012periodic,heemels2013model,borgers2018periodic,wang2019periodic} and STC \cite{mazo2009self,anta2010sample,wan2020dynamic,yi2018dynamic} of \textit{ODE systems} have surfaced. However, their use in controlling \textit{PDE plants} remains limited, with only a few papers addressing infinite-dimensional systems \cite{wakaiki2019stability,wakaiki2020event}. Utilizing \textit{semigroup theory,} \cite{wakaiki2020event} provides a \textit{full-state feedback} PETC for infinite dimensional systems with unbounded control operators and point actuation whereas \cite{wakaiki2019stability} provides a \textit{full-state feedback} STC for infinite dimensional systems with \textit{bounded control operators and spatially distributed actuation.}

This contribution introduces the first \textit{observer-based} PETC and STC for \textit{boundary control} of a class of parabolic PDEs using PDE backstepping approach\footnotemark. Specifically, we present observer-based PETC designs when the boundary sensing and actuation are either collocated or anti-collocated and an observer-based STC design with collocated boundary sensing and actuation. Our designs are far from trivial and encompass all possible configurations of boundary sensing and actuation but anti-collocated sensing and actuation under STC.
\footnotetext{Full-state feedback PETC and STC with global exponential convergence results were presented in \cite{rathnayake2023periodic} and \cite{rathnayake2023self}, respectively.}

The PETC results from a careful redesign of the continuous-time triggering function used in the CETC \cite{rathnayake2021observer,rathnayake2022sampled} to allow for periodic evaluation only. CETC approaches endowed with MDTs, such as those described in \cite{rathnayake2021observer,rathnayake2022sampled,espitia2020observer,espitia2019event,katz2020boundary,wang2022event,rathnayake2022observer}, might seem operable in a PETC configuration by permitting only periodic evaluations of the CETC triggering function, with a period that is less than or equal to the MDT. Indeed, following a control update triggered by an event, it is not necessary to check the CETC triggering function until a period equivalent to the MDT has elapsed, because during this time, the function continues to satisfy the required condition. However, once this period has passed, the triggering function could violate the required condition at any moment, indicating that continuous evaluation of the triggering function is \textit{absolutely} necessary after a period equivalent to the MDT following an event to ensure that the condition is continuously met. This implies that the simplistic approach of periodically checking the CETC triggering function at intervals less than or equal to the MDT does not suffice to guarantee the continuous satisfaction of the required condition. Thus, novel triggering functions designed for periodic evaluation are necessary. 

We derive a novel periodic event triggering function requiring only periodic evaluations by finding an upper bound of the underlying continuous-time triggering function between two consecutive periodic evaluations. Subsequently, an explicit upper-bound of the allowable sampling period for periodic evaluation of the triggering function is obtained. Since the triggering function is evaluated periodically, and the control input is updated only when the function satisfies a certain condition upon evaluation, Zeno behavior is inherently absent. Moreover, it is rigorously proven that the closed-loop system well-posedness and convergence under the CETC are preserved under the PETC. Specifically, the closed-loop signals under both CETC and PETC globally exponentially converges to zero in the spatial $L^2$ norm at the same rate. 

The proposed STC consists of a uniformly and positively lower-bounded function that accepts several inputs involving the observer states, which, when evaluated at an event time, outputs the waiting time until the next event. The design of the positive function requires upper and lower bounds of constituent variables  of the triggering function of the CETC. Since the function is uniformly and positively lower-bounded, the closed-loop system is Zeno-free by design. Moreover, the well-posedness of the closed-loop system under the STC is established. It is also proven that the closed-loop system under the STC exponentially converges to zero in the spatial $L^2$ norm locally at the same rate as its CETC counterpart.

\textit{Notation:} By $C^{0}(A;\Omega)$, we denote the class of continuous functions on $A\subseteq\mathbb{R}^{n}$, which takes values in $\Omega\subseteq\mathbb{R}$. By $C^{k}(A;\Omega)$, where $k\geq 1$, we denote the class of continuous functions on $A$, which takes values in $\Omega$ and has continuous derivatives of order $k$.  $L^{2}(0,1)$ denotes the equivalence class of Lebesgue measurable functions $f:[0,1]\rightarrow\mathbb{R}$ such that $\Vert f\Vert=\big(\int_{0}^{1}\vert f(x)\vert^{2}\big)^{1/2}<\infty$. $H^{1}(0,1)$ denotes the equivalence class of Lebesgue measurable functions $f:[0,1]\rightarrow\mathbb{R}$ such that $\int_{0}^{1} f^2(x)dx+\int_{0}^{1} f_x^2(x)dx<\infty$. Let $u:[0,1]\times\mathbb{R}_{+}\rightarrow\mathbb{R}$ be given. $u[t]$ denotes the profile of $u$ at certain $t\geq 0$, \textit{i.e.,} $\big(u[t]\big)(x)=u(x,t),$ for all $x\in [0,1]$. For an interval $J\subseteq\mathbb{R}_{+},$ the space $C^{0}\big(J;L^{2}(0,1)\big)$ is the space of continuous mappings $J\ni t\rightarrow u[t]\in L^{2}(0,1)$.

\section{Preliminaries and Continuous-time Event-triggered Control (CETC)}
Consider the following 1-D reaction-diffusion sampled-data boundary control system with constant coefficients:
\begin{align}\label{ctpe1} 
u_{t}(x,t)&=\varepsilon u_{xx}(x,t)+\lambda u(x,t),\text{ for }x\in(0,1),\\
\label{ctpe2}
\theta_1u_{x}(0,t)&=-\theta_2u(0,t),\\\label{ctpe3}
u_{x}(1,t)&=-qu(1,t)+U^{\omega}_j,
\end{align}
for all $t\in (t^\omega_j,t^\omega_{j+1}),j\in\mathbb{N}$, where $\theta_1\theta_2=0,\ \theta_1+\theta_2=1$, $``\omega"\in\{``c",``p",``s"\},$ and $t^\omega_0=0$. The sets $\{t^c_{j}\}_{j\in\mathbb{N}}$, $\{t^p_{j}\}_{j\in\mathbb{N}}$, and $\{t^s_{j}\}_{j\in\mathbb{N}}$ are event sequences from continuous-time event-triggering, periodic event-triggering, and self-triggering mechanisms. The initial condition is $u[0]\in L^{2}(0,1)$, and the parameters $\varepsilon,\lambda,q$ are all positive. The inputs $U_j^c,U_j^p,$ and $U_j^s$ are CETC, PETC, and STC inputs, respectively, held constant for $t\in[t_j^\omega,t^\omega_{j+1}),j\in\mathbb{N}$. Note that $\theta_1$ and $\theta_2$ are either $0$ or $1$ and $\theta_1\neq \theta_2$. The case $\theta_1=1,\theta_2=0$ leads to Neumann boundary condition at $x=0$ whereas $\theta_1=0,\theta_2=1$ leads to Dirichlet boundary condition at $x=0$.

In \cite{rathnayake2021observer} and \cite{rathnayake2022sampled}, the authors develop observers for the system \eqref{ctpe1}-\eqref{ctpe3} using boundary measurements. The former addresses anti-collocated boundary sensing and actuation with $u(0,t)$ as the measurement, while the latter focuses on collocated boundary sensing and actuation with $u(1,t)$ as the measurement. These designs are presented below:

\begin{align}\label{ctoe1kl}
\begin{split}
\hat{u}_{t}(x,t)=&\varepsilon \hat{u}_{xx}(x,t)+\lambda \hat{u}(x,t)\\&+p_1(x)\big(\theta_1\tilde{u}(0,t)+\theta_2\tilde{u}(1,t)\big),\text{ for }x\in(0,1),
\end{split}
\\\label{ctoe2kl}
\theta_1\hat{u}_x(0,t)&=-\theta_2\hat{u}(0,t)+\theta_1p_{10}\tilde{u}(0,t),
\\\label{ctoe3kl}
\hat{u}_{x}(1,t)&=-q\hat{u}(1,t)+U^\omega_j+\theta_2p_{10}\tilde{u}(1,t),
\end{align}
for all $t\in (t^\omega_j,t^\omega_{j+1}),j\in\mathbb{N}$, where $\hat{u}[0]\in L^{2}(0,1)$ and
\begin{equation}\label{aqwr}
    \tilde{u}(x,t):=u(x,t)-\hat{u}(x,t),
\end{equation}
is the observer error. Here, the case $\theta_1=1,\theta_2=0$ results in anti-collocated sensing and actuation, and the case $\theta_1=0,\theta_2=1$ results in collocated sensing and actuation. The terms $p_1(x)$ and $p_{10}$ are observer gains determined using the PDE backstepping technique equipped with the Volterra transformation: 
\begin{equation}\label{qppu}
    \tilde{u}(x,t)=\tilde{w}(x,t)-\int_{\theta_2x}^{\theta_1x+\theta_2}P(x,y)\tilde{w}(y,t)dy,
\end{equation}
with its inverse:
\begin{equation}\label{qmpo}
    \tilde{w}(x,t)=\tilde{u}(x,t)+\int_{\theta_2x}^{\theta_1x+\theta_2}Q(x,y)\tilde{u}(y,t)dy,
\end{equation}
for $0\leq \theta_2x+\theta_1y\leq \theta_1x+\theta_2y\leq 1$. Details on the bounded observer gains $p_1(x)$ and $p_{10}$ as well as the bounded gain kernels $P(x,y)$ and $Q(x,y)$ are found in 
\cite{rathnayake2021observer} and \cite{rathnayake2022sampled}.

The well-posedness of the closed-loop system \eqref{ctpe1}-\eqref{aqwr} with piecewise constant inputs between two sampling instants is provided in the following proposition. 

\begin{prop}[\cite{rathnayake2021observer}]\label{cor1}
For every $u[t_{j}^\omega],\hat{u}[t_{j}^\omega]\in L^{2}(0,1)$, there exist unique solutions $u,\hat{u}:[t_j^\omega,t_{j+1}^\omega]\times[0,1]\rightarrow\mathbb{R}$ between two time instants $t_{j}^\omega$ and $t_{j+1}^\omega$ such that $u,\hat{u}\in C^{0}([t_{j}^\omega,t_{j+1}^\omega];L^{2}(0,1))\cap C^{1}((t_{j}^\omega,t_{j+1}^\omega)\times [0,1])$ with $u[t],\hat{u}[t]\in C^{2}([0,1])$ which satisfy \eqref{ctpe2},\eqref{ctpe3},\eqref{ctoe2kl},\eqref{ctoe3kl} for $t\in (t_{j}^\omega,t_{j+1}^\omega]$ and \eqref{ctpe1}, \eqref{ctoe1kl} for $t\in (t_{j}^\omega,t_{j+1}^\omega]$, $x\in(0,1)$.    
\end{prop}

\begin{ass}\label{ass1} The parameters $q,\lambda,$ and $\varepsilon$ satisfy the following relation:
\begin{equation}
    q>\frac{\lambda}{2\varepsilon}+\frac{\theta_1}{2}.
\end{equation}
\end{ass}


\begin{rmk}\label{rem1}\rm
Assumption \ref{ass1} is required to avoid the use of the signal $u(1, t)$ in the nominal control law for which it is impossible to obtain a useful bound on its rate of change. Furthermore, It is worth mentioning that an eigenfunction expansion of the solution of \eqref{ctpe1}-\eqref{ctpe3} with $U_j^\omega=0$ (zero input) shows that the system is unstable when $\lambda>\varepsilon\pi^2/4^{\theta_1}$, no matter what $q>0$ (see Remark 1 in \cite{rathnayake2021observer} and \cite{rathnayake2022sampled}). \hfill $\square$
\end{rmk}


In \cite{rathnayake2021observer} and \cite{rathnayake2022sampled}, the authors propose the following sampled-data boundary control law to be used in conjunction with event-triggering:
\begin{equation}\label{etcl}
U_{j}^\omega:=\int_{0}^{1}k(y)\hat{u}(y,t^\omega_{j})dy,
\end{equation}
for $t\in [t_j^\omega,t_{j+1}^\omega),j\in\mathbb{N}$ where $k(y)$ is the control gain found using the PDE backstepping technique equipped with the Volterra transformation: 
\begin{equation}\label{ctbt}
\hat{w}(x,t)=\hat{u}(x,t)-\int_{0}^{x}K(x,y)\hat{u}(y,t)dy,
\end{equation}
with its inverse:
\begin{equation}\label{ink}
\hat{u}(x,t)=\hat{w}(x,t)+\int_{0}^{x}L(x,y)\hat{w}(y,t)dy,
\end{equation}
for $0\leq y\leq x\leq 1$. For further details on the bounded control gain $k(x)$ as well as the bounded gain kernels $K(x,y)$ and $L(x,y)$, the readers are referred to \cite{rathnayake2021observer} and \cite{rathnayake2022sampled}. 

The difference between the sampled-data and the continuous-time control input, termed the input holding error, is defined by
\begin{equation}\label{dt}
d(t):=\int_{0}^{1}k(y)\big(\hat{u}(y,t_j^\omega)-\hat{u}(y,t)\big)dy,
\end{equation}
where $t\in[t_j^\omega,t_{j+1}^{\omega})$ and $j\in\mathbb{N}$.

The authors of \cite{rathnayake2021observer,rathnayake2022sampled} present a continuous-time event-triggering mechanism to determine the set of event times $\{t_j^{c}\}_{j\in\mathbb{N}}$ using $d(t)$ and a dynamic variable $m(t)$ via the following rule:
\begin{equation}\label{tnxt}
 t_{j+1}^c=\inf\big\{t\in\mathbb{R}_+\vert t>t_j^c,\Gamma^c(t)>0,j\in\mathbb{N}\big\},
\end{equation}
with $t^c_0=0$. The function $\Gamma^c(t)$ is defined as 
\begin{equation}\label{gmmt}
  \Gamma^c(t):=d^2(t)-\gamma m(t),
\end{equation}
where $\gamma>0$ is an event-trigger design parameter. The variable $m(t)$ evolves according to the ODE
\begin{equation}\label{obetbc3m}\begin{split}
\dot{m}(t)=&-\eta m(t)-\rho d^{2}(t)+\beta_{1}\Vert \hat{u}[t]\Vert^{2}+\beta_{2}\hat{u}^2(1,t)\\&+\theta_1\beta_3\tilde{u}^2(0,t)+\theta_2\beta_3\tilde{u}^2(1,t),\end{split}
\end{equation}
valid for all $t\in(t_{j}^c,t_{j+1}^c),j\in\mathbb{N}$ with $m(t_{0}^c)=m(0)>0$ and $m(t_{j}^{c-})=m(t^c_{j})=m(t_{j}^{c+})$, and $\eta,\rho,\beta_1,\beta_2,\beta_3>0$ are event-trigger parameters.

\begin{ass}[Event-trigger parameter selection]\label{ass2}
     The parameters $\gamma,\eta>0$ are design parameters, and $\beta_1,\beta_2,\beta_3>0$ are chosen such that
\begin{equation}\label{betasj}
\beta_{1}=\frac{\alpha_{1}}{\gamma(1-\sigma)},\hspace{5pt}\beta_{2}=\frac{\alpha_{2}}{\gamma(1-\sigma)},\hspace{5pt} \beta_{3}=\frac{\alpha_{3}}{\gamma(1-\sigma)},
\end{equation}
where $\sigma\in(0,1)$ and 
\begin{align}
\begin{split}
\label{al1j}\alpha_{1}&=4\int_{0}^{1}\Big(\varepsilon k''(y)+\varepsilon k(1)k(y)+\lambda k(y)\Big)^{2}dy,
\end{split}\\
\label{al2j}
\alpha_{2}&=4\big(\varepsilon qk(1)+\varepsilon k'(1)\big)^{2},\\\label{alj3}
\alpha_3 &= 4\Big(\frac{\lambda \big(\theta_1k(0)
+\theta_2k(1)\big)}{2}+\int_{0}^{1}k(y)p_1(y)dy\Big)^2,
\end{align}
Subject to Assumption \ref{ass1}, the parameter $\rho>0$ is chosen as 
\begin{equation}\label{hhjk}
\rho=\frac{\varepsilon\kappa_1 B}{2},
\end{equation}
for $B,\kappa_1>0$ chosen such that $B\Big(\varepsilon\min\Big\{q-\frac{\lambda}{2\varepsilon}-\frac{\theta_1}{2},\frac{1}{2}\Big\}-\frac{\varepsilon}{2\kappa_1}-\frac{\lambda\big(5\theta_1+2\theta_2\big)}{8\kappa_2}-\frac{\Vert g\Vert^2}{\kappa_3}\Big)-2\beta_1\tilde{L}^2-2\beta_{2}-4\beta_{2}\check{L}^2>0,$ for some $\kappa_2,\kappa_3>0$, where $g(x)=p_{1}(x)-\frac{\theta_1\lambda}{2}K(x,0)-\int_{0}^{x}K(x,y)p_{1}(y)dy$, $\tilde{L} = 1+\big(\int_{0}^{1}\int_{0}^xL^2(x,y)dydx\big)^{1/2}$ and $\check{L}=\big(\int_{0}^1L^2(1,y)dy\big)^{1/2}$ with $K(x,y)$ and $L(x,y)$ being the gain kernels of the backstepping transformations \eqref{ctbt} and \eqref{ink}. Note from Assumption \ref{ass1} that $q-\lambda/2\varepsilon-\theta_1/2>0$. 
\end{ass}

\begin{thme}[Results under CETC \cite{rathnayake2021observer,rathnayake2022sampled}]
\label{hhgb} Consider the CETC approach \eqref{etcl},\eqref{dt}-\eqref{obetbc3m} under Assumption \ref{ass1}, which generates a set of event-times $I^c=\{t_j^c\}_{j\in\mathbb{N}}$ with $t_0^c=0$. It holds that
\begin{equation}
    \Gamma^c(t)\leq 0 \text{ for } \text{ all }t\in \big[0,\sup(I^c)\big),
\end{equation}
Consequently, given appropriate choices for the event-trigger parameters $\gamma$,$\eta$,$\beta_1$,$\beta_2$,$\beta_3$,$\rho>0$, the followings hold:
\begin{enumerate}
    \item[R1:] The set of event-times $I^c$ generates an increasing sequence for any $\eta,\gamma,\rho>0$ and $\beta_1,\beta_2,\beta_3>0$ satisfying \eqref{betasj}. Specifically, it holds that $t^c_{j+1}-t^c_j\geq \tau>0,j\in\mathbb{N}$ where
    \begin{equation}\label{mdt}
\tau=\frac{1}{a}\ln \Big ( 1+\frac{\sigma a}{(1-\sigma)(a+\gamma\rho)}\Big).
\end{equation}
Here $\sigma\in(0,1)$ appears in the relation \eqref{betasj}, and \begin{equation}\label{hhnm}
    a = 1+\rho_1+\eta>0,
\end{equation}
where 
\begin{equation}\label{sopw}
    \rho_1=4\varepsilon^{2}k^{2}(1).
\end{equation}
As $j\rightarrow \infty$, it follows that $t_j^c\rightarrow \infty$, thereby excluding Zeno behavior.
    \item [R2:] For every $u[0],\hat{u}[0]\in L^{2}(0,1)$, there exist unique solutions $u,\hat{u}:\mathbb{R}_+\times[0,1]\rightarrow\mathbb{R}$ such that $u,\hat{u}\in C^0(\mathbb{R}_+;L^2(0,1)\cap C^1(J^c\times [0,1])$ with $u[t],\hat{u}[t]\in C^{2}([0,1])$ which satisfy \eqref{ctpe2},\eqref{ctpe3},\eqref{ctoe2kl},\eqref{ctoe3kl} for all $t>0$ and \eqref{ctpe1}, \eqref{ctoe1kl} for all $t>0,x\in(0,1),$ where $J^c=\mathbb{R_{+}}\text{\textbackslash}I^c$.
\item[R3:] The dynamic variable $m(t)$ governed by \eqref{obetbc3m} with $m(0)>0$ satisfies $m(t)>0$ for all $t>0$.

\item[R4:] As a result of R1-R3 and under Assumption \ref{ass2}, the closed-loop system \eqref{ctpe1}-\eqref{aqwr} globally exponentially converges to zero in the spatial $L^2$ norm satisfying
\begin{equation}\label{rrt}
    \Vert u[t]\Vert+\Vert\hat{u}[t]\Vert\leq M e^{-\frac{b^{*}}{2}t}\sqrt{\Vert u[0]\Vert^2+\Vert \hat{u}[0]\Vert^2+m(0)},
\end{equation}
for all $t>0$ and for some $M,b^*>0$.
\end{enumerate}
\end{thme}

\section{Periodic Event-triggered Control (PETC) and Self-triggered Control (STC)}

\subsection{Periodic Event-triggered Control (PETC)}
In this subsection, we propose a PETC approach for the system described by equations \eqref{ctpe1}-\eqref{aqwr} subject to Assumption \ref{ass1}. This approach is applicable under both anti-collocated ($\theta_1=1,\theta_2=0$) and collocated ($\theta_1=0,\theta_2=1$) sensing and actuation configurations. Our design draws from the CETC scheme in \eqref{etcl}, \eqref{dt}-\eqref{obetbc3m}. To implement this, we redesign the triggering function $\Gamma^c(t)$, as specified in \eqref{gmmt}, into a new triggering function $\Gamma^p(t)$ which facilitates periodic evaluations. Furthermore, we determine a maximum allowable sampling period $h>0$ for the periodic event-trigger. The proposed periodic event-triggering mechanism determines the set of event-times $\{t_j^p\}_{j\in\mathbb{N}}$ via the following rule:
\begin{equation}\label{petg}
    \begin{split}
    t_{j+1}^p=\inf\big\{t\in\mathbb{R}_+\vert& t>t_j^p,\Gamma^p(t)>0, t=nh,\\&\qquad\qquad h>0,n\in\mathbb{N}\big\},
    \end{split}
\end{equation}
with $t^p_0=0$. Here, $h$ is the sampling period selected as
\begin{equation}\label{hjgf}
    0< h\leq \tau,
\end{equation}
where $\tau$ is given by \eqref{mdt} and $\Gamma^p(t)$ is given by
\begin{equation}\label{tildeGqwref}
\begin{split}
\Gamma^p(t )=&
(a+\gamma\rho)e^{ah}d^{2}(t )-\gamma\rho d^2(t )-\gamma a m(t).
\end{split}
\end{equation}
Here, $d(t)$ is given by \eqref{dt} for $t\in[t_j^p,t_{j+1}^p),j\in\mathbb{N}$, $m(t)$ satisfies \eqref{obetbc3m} for $t\in(t_j^p,t_{j+1}^p),j\in\mathbb{N}$, and $a$ is given by \eqref{hhnm}.

Note that, under the continuous-time event-trigger \eqref{tnxt}-\eqref{obetbc3m}, the triggering function $\Gamma^c(t)$ needs to be checked continuously to detect events. In contrast, with the periodic event-trigger \eqref{petg}-\eqref{tildeGqwref}, the triggering function $\Gamma^p(t)$ requires only periodic evaluations for event detection. 

 \begin{thme}[Results under PETC]\label{hhgbnml} Consider the PETC approach \eqref{etcl},\eqref{petg}-\eqref{tildeGqwref} under Assumption \ref{ass1}, which generates an increasing set of event-times $I^p=\{t_j^p\}_{j\in\mathbb{N}}$ with $t_0^p=0$. For every $u[0],\hat{u}[0]\in L^{2}(0,1)$, there exist unique solutions $u,\hat{u}:\mathbb{R}_+\times[0,1]\rightarrow\mathbb{R}$ such that $u,\hat{u}\in C^0(\mathbb{R}_+;L^2(0,1)\cap C^1(J^p\times [0,1])$ with $u[t],\hat{u}[t]\in C^{2}([0,1])$ which satisfy \eqref{ctpe2},\eqref{ctpe3},\eqref{ctoe2kl},\eqref{ctoe3kl} for all $t>0$ and \eqref{ctpe1}, \eqref{ctoe1kl} for all $t>0,x\in(0,1),$ where $J^p=\mathbb{R_{+}}\text{\textbackslash}I^p$. Given appropriate choices for the event-trigger parameters $\gamma$,$\eta$,$\beta_1$,$\beta_2,$$\beta_3$,$\rho>0$, the followings hold: 
\begin{enumerate}
    \item[R1:] For any $\eta,\gamma,\rho>0$ and $\beta_1,\beta_2,\beta_3>0$ satisfying \eqref{betasj}, the function $\Gamma^c(t)$ given by \eqref{gmmt}  satisfies $\Gamma^c(t)\leq 0$ for all $t>0$ along the solution of \eqref{ctpe1}-\eqref{aqwr},\eqref{etcl},\eqref{dt},\eqref{obetbc3m},\eqref{petg}-\eqref{tildeGqwref}.   
    \item[R2:] The dynamic variable $m(t)$ governed by \eqref{obetbc3m} with $m(0)>0$ satisfies $m(t)>0$ for all $t>0$ along the solution of \eqref{ctpe1}-\eqref{aqwr},\eqref{etcl},\eqref{dt},\eqref{petg}-\eqref{tildeGqwref}.
    \item[R3:] Under Assumption \ref{ass2}, the closed-loop system \eqref{ctpe1}-\eqref{aqwr} globally exponentially converges to zero in the spatial $L^2$ norm satisfying \eqref{rrt}.
\end{enumerate}
\end{thme}
\noindent The complete proof is provided in the Appendix-A.

\subsection{Self-triggered Control (STC)}
In this subsection, we propose an STC approach for the system described by equations \eqref{ctpe1}-\eqref{aqwr} under collocated sensing and actuation configuration ($\theta_1=0,\theta_2=1$)
 and subject to Assumption \ref{ass1}. Furthermore, we make the following assumption on initial data.
\begin{ass}\label{ddcf}
    The initial conditions of the plant \eqref{ctpe1}-\eqref{ctpe3} and the observer \eqref{ctoe1kl}-\eqref{ctoe3kl} with $\theta_1=0,\theta_2=1$ satisfy $u[0],\hat{u}[0]\in H^1(0,1)$. Further, for some known constants $\Psi_1,\Psi_2>0,$ it holds that
    \begin{equation}\label{hjg1}
        \Vert u[0]\Vert\leq \Psi_1, \text{ and }\Vert u_x[0]\Vert\leq \Psi_2.
    \end{equation}
\end{ass}

Our design draws from the CETC scheme in \eqref{etcl}, \eqref{dt}-\eqref{obetbc3m}. We propose a function $G(\cdot,\cdot)$ that is uniformly positive and lower-bounded, and which depends on the observer states. When this function is evaluated at the current event time, it yields the waiting time until the subsequent event. The proposed self-triggering mechanism determines the sequence of event times $\{t_j^s\}_{j\in\mathbb{N}}$ according to the following rule:
\begin{equation}\label{xxvb}
t^s_{j+1}=t^s_j+G\big(\Vert \hat{u}[t^s_j]\Vert,m(t^s_j)\big),
\end{equation}
 with $t^s_0=0$ where $G(\cdot,\cdot)>0$ is a  uniformly and positively lower-bounded function
  \begin{equation}\label{cvbs}
   \begin{split}
       &G(\Vert \hat{u}[t^s_j]\Vert,m(t^s_j))\\&:=\max\bigg\{\tau,\frac{1}{2\varrho+\eta}\ln\Big(\frac{\gamma m(t^s_j)+\frac{\gamma\rho H(t^s_j)}{2\varrho+\eta}}{H(t^s_j)+\frac{\gamma\rho H(t^s_j)}{2\varrho+\eta}}\Big)\bigg\},
    \end{split}
\end{equation}
 In \eqref{cvbs}, $H(t)$ is given by
   \begin{equation}\label{xxcvb}
    H(t)=2\Vert k\Vert^2\Big(2\Vert\hat{u}[t]\Vert^2+\frac{\varepsilon^2\Vert k\Vert^2}{\lambda\varrho}\Vert\hat{u}[t]\Vert^2+\frac{(\Psi_0^*)^2e^{-2\sigma^*t}}{\varrho}\Big),
\end{equation}
where
\begin{equation}\label{qqqwwe}
    \varrho = \lambda+\frac{\Vert p_1\Vert^2}{2},\text{ }\sigma^* \in \Big(0,\frac{\varepsilon\pi^2}{4}\Big],
\end{equation}
and $\Psi_0^*$ is given by
\begin{equation}\label{ytre}
\begin{split}    \Psi_0^*&=\frac{1}{\sqrt{2}}\bigg(\Big((M_1+1)\Omega_1+\Omega_2\Big)\Psi_1+\Psi_2\\&\qquad\qquad+\Big((M_1+1)\Omega_1+\Omega_2\Big)\Vert\hat{u}[0]\Vert+\Vert\hat{u}_x[0]\Vert\bigg)\\&\quad\times
\sqrt{\Big(\frac{\varepsilon^2p_{10}^2}{\lambda}+\frac{1}{2}\Big)}.
\end{split}
\end{equation}
 In \eqref{ytre}, $\Psi_1,\Psi_2>0$ are the known bounds of $\Vert u[0]\Vert$ and $\Vert u_x[0]\Vert$, respectively, as stated in Assumption \ref{ddcf} and $M_1,\Omega_1,\Omega_2$ are given by 
 \begin{equation}\label{phmk1}
     M_{1}=2q+\frac{2\varepsilon q^{2}}{\varepsilon\pi^{2}/4+2\varepsilon q-\sigma^*},
 \end{equation}
 \begin{equation}\label{zzzs1}
\Omega_1=1+\sqrt{\int_{0}^{1}\int_{x}^1Q^2(x,y)dydx},
 \end{equation}
 and 
\begin{equation}\label{dsfg}
    \Omega_2=\max_{x\in[0,1]}\vert Q(x,x)\vert+\sqrt{\int_{0}^1\int_{x}^1 Q_x^2(x,y)dydx},
\end{equation}
where $Q(x,y)$ is the gain kernel of the inverse backstepping transformation \eqref{qmpo}. Referring to \cite{rathnayake2022sampled}, one can show that $\max_{x\in[0,1]}\vert Q(x,x)\vert=\lambda/2\varepsilon$. 

\begin{thme}[Results under STC]\label{hhgbsd} Consider the STC approach \eqref{etcl},\eqref{xxvb}-\eqref{dsfg} under Assumptions \ref{ass1} and \ref{ddcf}, which generates an increasing set of event-times $I^s=\{t_j^s\}_{j\in\mathbb{N}}$ with $t_0^s=0$. For every $u[0],\hat{u}[0]\in L^{2}(0,1)$, there exist unique solutions $u,\hat{u}:\mathbb{R}_+\times[0,1]\rightarrow\mathbb{R}$ such that $u,\hat{u}\in C^0(\mathbb{R}_+;L^2(0,1)\cap C^1(J^s\times [0,1])$ with $u[t],\hat{u}[t]\in C^{2}([0,1])$ which satisfy \eqref{ctpe2},\eqref{ctpe3},\eqref{ctoe2kl},\eqref{ctoe3kl} for all $t>0$ and \eqref{ctpe1}, \eqref{ctoe1kl} for all $t>0,x\in(0,1),$ where $J^s=\mathbb{R_{+}}\text{\textbackslash}I^s$. Given appropriate choices for the event-trigger parameters $\gamma$,$\eta$,$\beta_1$,$\beta_2$,$\beta_3$,$\rho>0$, the followings hold: 
\begin{enumerate} 
   \item[R1:] For any $\eta,\gamma,\rho>0$ and $\beta_1,\beta_2,\beta_3>0$ satisfying \eqref{betasj}, the function $\Gamma^c(t)$ given by \eqref{gmmt}  satisfies $\Gamma^c(t)\leq 0$ for all $t>0$ along the solution of \eqref{ctpe1}-\eqref{aqwr},\eqref{etcl},\eqref{dt},\eqref{obetbc3m},\eqref{xxvb}-\eqref{dsfg}.   
    \item[R2:] The dynamic variable $m(t)$ governed by \eqref{obetbc3m} with $m(0)>0$ satisfies $m(t)>0$ for all $t>0$ along the solution of \eqref{ctpe1}-\eqref{aqwr},\eqref{etcl},\eqref{dt},\eqref{xxvb}-\eqref{dsfg}.
    \item[R3:] Under Assumption \ref{ass2}, the closed-loop system \eqref{ctpe1}-\eqref{aqwr} locally exponentially converges to zero in the spatial $L^2$ norm satisfying \eqref{rrt}.
\end{enumerate}
\end{thme}
\noindent The complete proof is provided in the Appendix-B.
\begin{rmk}\rm\label{ici}
When $u[0]$ and $\hat{u}[0]$ are in $H^{1}(0,1)$, they are also in $L^{2}(0,1)$. As a result, the well-posedness of the closed-loop system \eqref{ctpe1}-\eqref{ctoe3kl}, in the context of Theorem \ref{hhgbsd}, directly follows from Proposition \ref{cor1}. The solution is constructed iteratively between consecutive event times.
\end{rmk}

\section{Numerical Simulations}
\begin{figure}
\centering
\includegraphics[scale=.075]{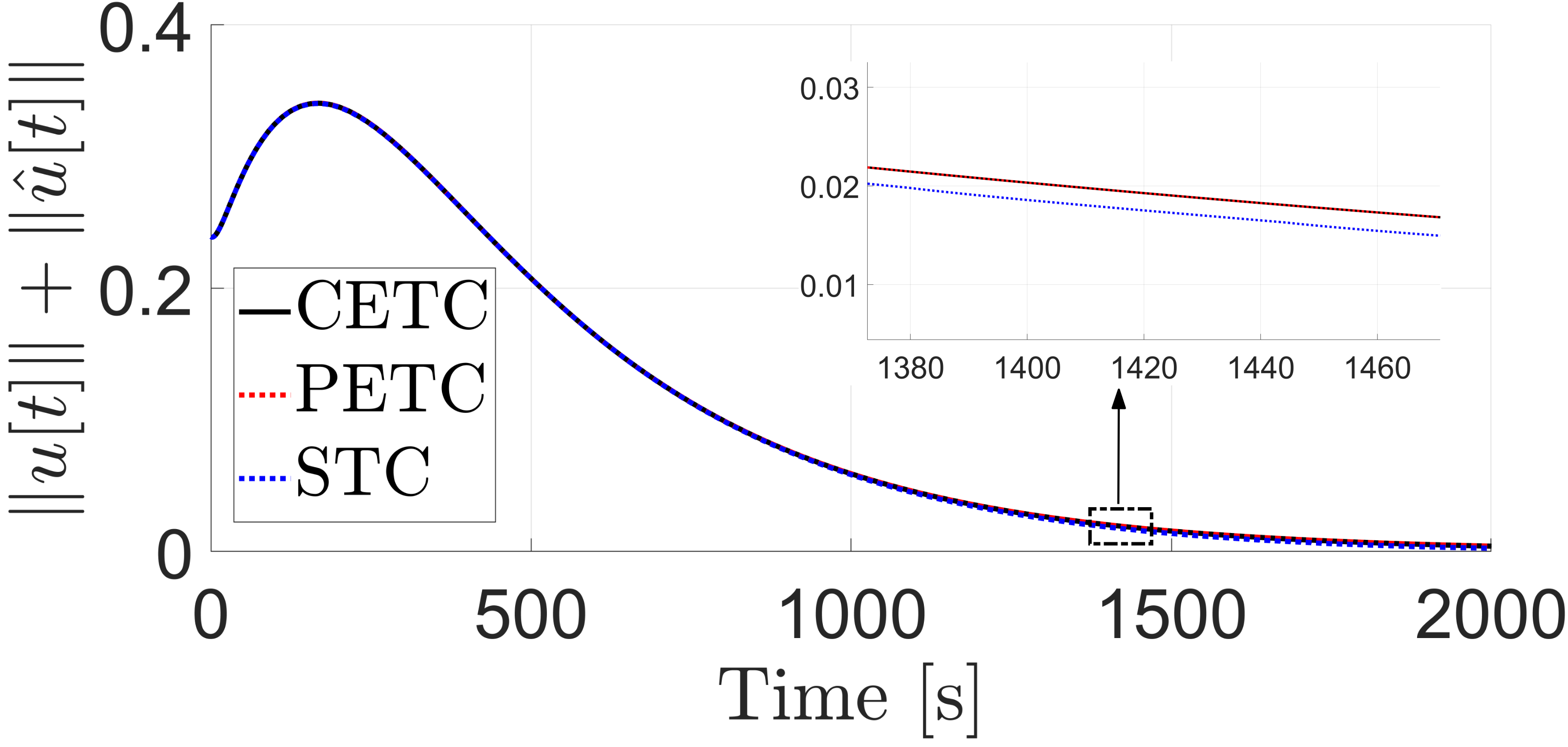}
\caption{Evolution of $\Vert u[t]\Vert+\Vert\hat{u}[t]\Vert$}.
\label{cv51}
\end{figure}

\begin{figure}
\centering
\includegraphics[scale=.075]{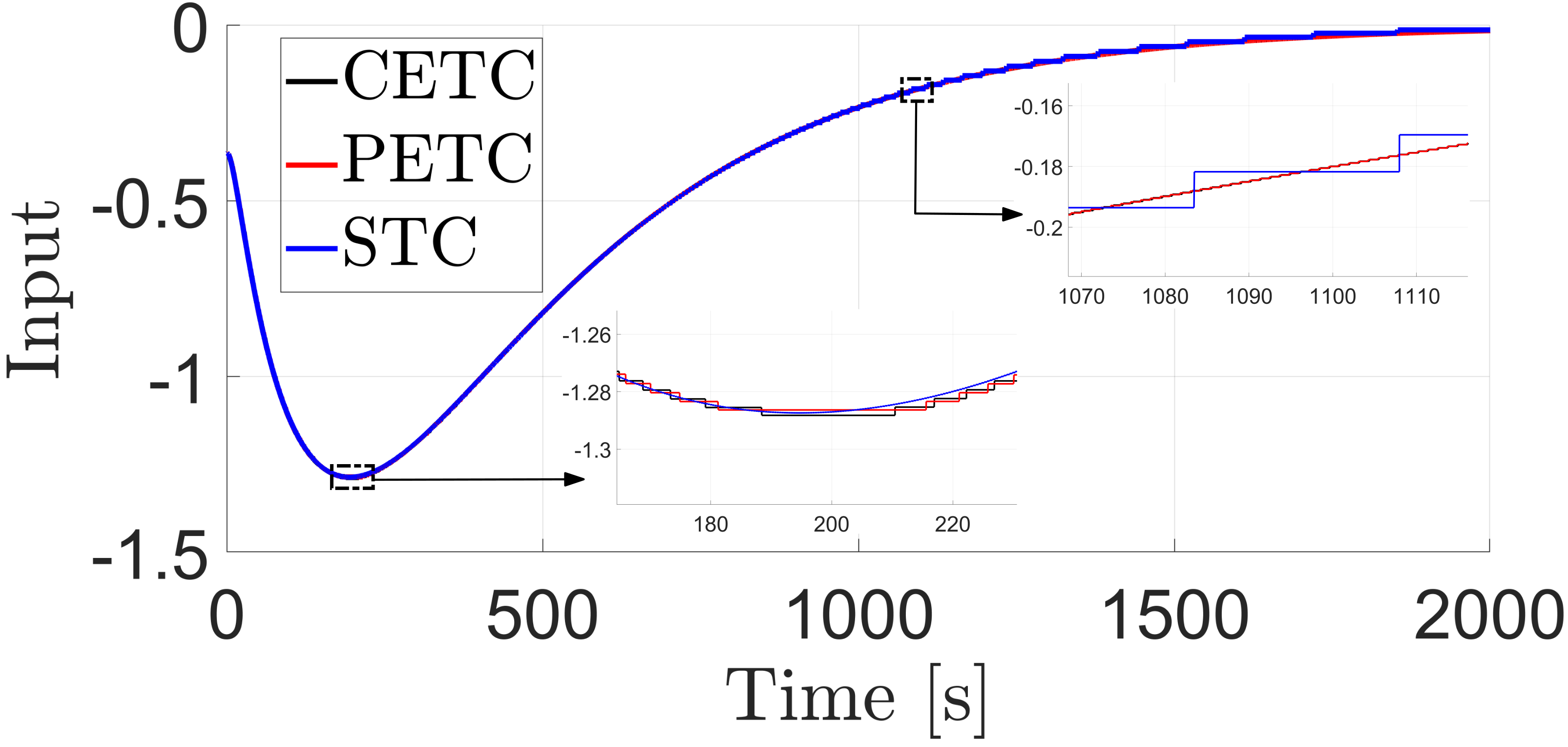}
\caption{Boundary control inputs.}
\label{cv52}
\end{figure}

\begin{figure}
\centering
\includegraphics[scale=.09]{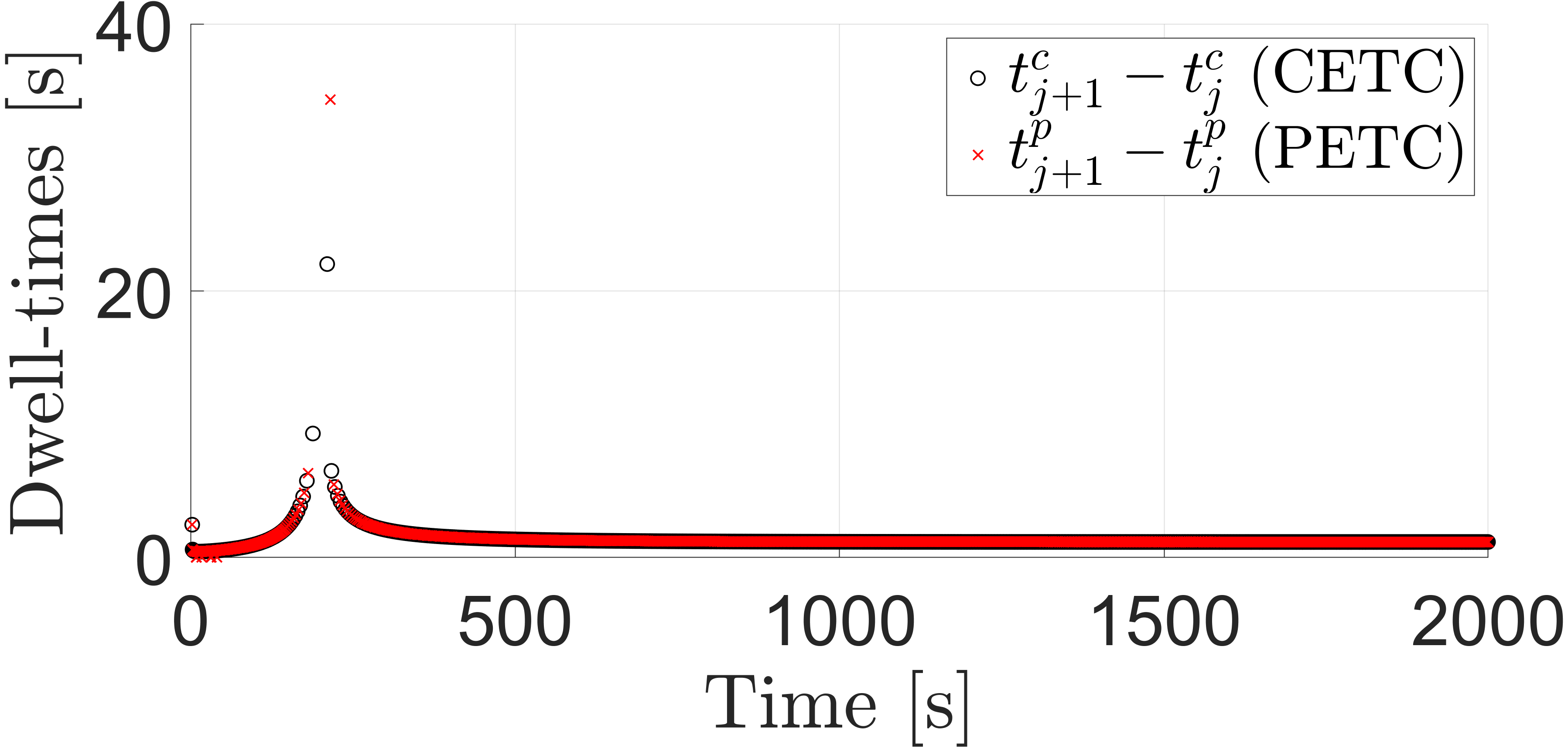}
\caption{Dwell-times under CETC and PETC.}
\label{cv53}
\end{figure}

\begin{figure}
\centering
\includegraphics[scale=.09]{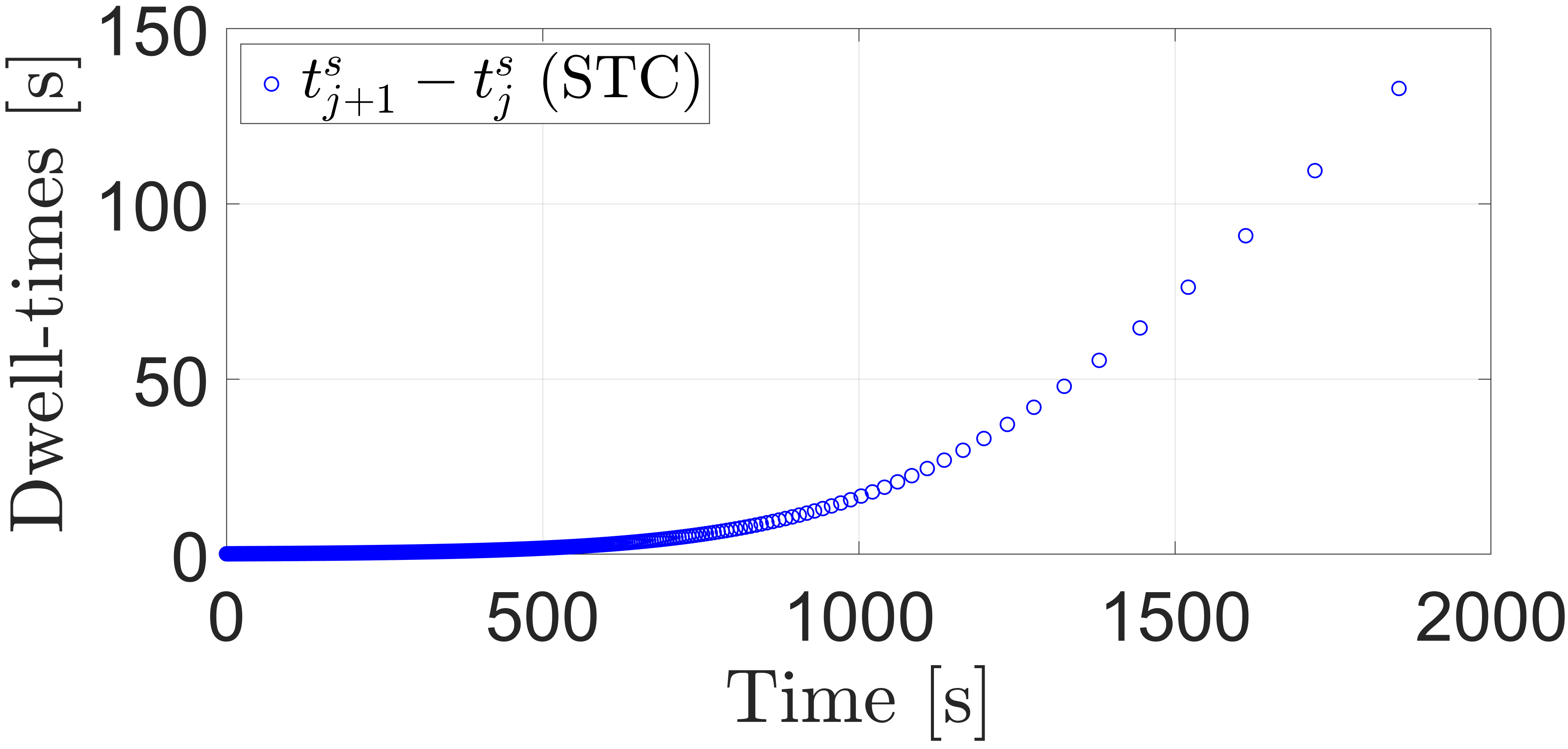}
\caption{Dwell-times under STC.}
\label{cv54}
\end{figure}

We consider an open loop unstable reaction-diffusion PDE with $\varepsilon=0.001,\lambda=0.01,q=5.1,\theta_1=0,\theta_2=1$ and the initial conditions $u[0]=5x^2(x-1)^{2}$ and $\hat{u}[0]=x^2(x-1)^{2}$. For numerical simulations, both the plant and the observer are discretized with a uniform step size of $\Delta x=0.005$ for the space variable. The discretization with respect to time is done using the implicit Euler scheme with step size $\Delta t=0.001s$. The parameters for the CETC and PETC are chosen as follows: $m(0)=10^{-4},\gamma=1,\eta=1$ and  $\sigma=0.9$. We compute using \eqref{al1j}-\eqref{alj3} that $\alpha_{1}= 0.021;\alpha_2= 0.0131;\alpha_3=0.7971$. Therefore, from \eqref{betasj}, we obtain $\beta_{1}= 0.2095;\beta_{2}= 0.1309;\beta_{3}=  7.9706$. Let us choose $\kappa_1=25$ and $B=7.7304\times 10^4$, and then, from \eqref{hhjk}, we obtain $\rho=966.3$. The parameters for the STC are chosen as follows: $m(0)=10^{-4},\gamma=10^{12},\eta=10^{-6}, \sigma=0.9,\Psi_1=0.1992,\Psi_2=0.6901,$ and $\sigma^*=\varepsilon\pi^2/4$. Therefore, from \eqref{betasj}, we obtain $\beta_{1}=2.095\times 10^{-13};\beta_{2}= 1.31\times 10^{-13};\beta_{3}=   7.9706 \times 10^{-12}$. Let us choose $\kappa_1=25$ and $B= 7.7304\times 10^{-8}$, and then, from \eqref{hhjk}, we obtain $\rho=9.6630\times 10^{-10}$. For all the strategies, the MDT $\tau$ calculated using \eqref{mdt} is $0.009s$. For the PETC, we choose the sampling period for evaluating the triggering function as $h=0.009s$. 

Fig. \ref{cv51} shows the response of the closed-loop system under CETC, PETC, and STC. Fig. \ref{cv52} shows the corresponding control inputs. We can observe that the spatial $L^2$ norm of the closed-loop signals under CETC, PETC, and STC converge to zero at roughly similar rates, despite PETC and STC not requiring continuous evaluation of a triggering function, as opposed to CETC. In Figs. \ref{cv53} and \ref{cv54}, we can observe that both CETC and PETC trigger events at similar rates, whereas STC triggers more frequent events than CETC and PETC initially. However, as the closed-loop system approaches the equilibrium, the frequency of events under STC significantly becomes lower than those under CETC and PETC. 

\section{Conclusions}
In this paper, we have proposed novel observer-based PETC and STC strategies for a class of reaction-diffusion systems. Specifically, we have presented observer-based PETC designs when the sensing and actuation are either collocated or anti-collocated and  an observer-based STC design with collocated sensing and actuation. The key idea of the developed methods is the transformation of a class of continuous-time dynamic event-triggers which require continuous evaluation to periodic event-triggers and self-triggers. For the PETC, we have obtained an explicit upper-bound of the allowable sampling period of the periodic event-trigger. For the STC, we have designed a uniformly and positively lower-bounded function which, when evaluated at the time of an event, outputs the waiting time until the next event. The well-posedness of the closed-loop system under both PETC and STC has been proven for all cases. Further, we have proven that the global exponential convergence  to  zero in the spatial $L^2$ norm under the CETC is preserved under the proposed PETC. The STC ensures that the closed-loop system exponentially converges to zero in the spatial $L^2$ norm locally at a comparable rate to its CETC counterpart. The conducted numerical simulation has illustrated the validity of the theoretical developments.

\section*{Appendix}

\noindent \textbf{A. Proof of Theorem \ref{hhgbnml}}

 The well-posedness of the closed-loop system \eqref{ctpe1}-\eqref{ctoe3kl} under the PETC, in the sense of Theorem \ref{hhgbnml}, directly follows from Proposition \ref{cor1}. The solution is constructed iteratively between consecutive event times. To streamline the rest of the proof of Theorem \ref{hhgbnml}, we first present Lemmas \ref{lem2} and \ref{imptnt}.  

\begin{lem}\label{lem2} Consider the PETC approach \eqref{etcl},\eqref{petg}-\eqref{tildeGqwref} which generates an increasing set of event-times $I^p=\{t_j^p\}_{j\in\mathbb{N}}$ with $t^p_0=0$.  For $d(t)$ given by \eqref{dt}, it holds that\begin{equation}\label{ghm}
\begin{split}
\dot{d}^{2}(t)\leq &\rho_{1} d^{2}(t)+\alpha_{1}\Vert \hat{u}[t]\Vert^{2}+\alpha_{2}\hat{u}^2(1,t)+\theta_1\alpha_3\tilde{u}^2(0,t)\\&+\theta_2\alpha_3\tilde{u}^2(1,t),
\end{split}
\end{equation}
along the solution of  \eqref{ctpe1}-\eqref{aqwr} for all $t\in\big(nh,(n+1)h\big)$ and any $n\in\big[t^p_j/h,t^p_{j+1}/h\big)\cap\mathbb{N}$. Here $\alpha_1,\alpha_2,\alpha_3,\rho_1>0$ are given by \eqref{al1j}-\eqref{alj3},\eqref{sopw}, respectively.\end{lem}

The proof is very similar to that of Lemma 2 in \cite{rathnayake2021observer}, and hence omitted.

\begin{lem}\label{imptnt}  Consider the PETC approach \eqref{etcl},\eqref{petg}-\eqref{tildeGqwref} under Assumption \ref{ass1}, which generates an increasing set of event-times $\{t^p\}_{j\in\mathbb{N}}$ with $t^p_0=0$. For any $\eta,\gamma,\rho>0$ and $\beta_1,\beta_2,\beta_3>0$ satisfying \eqref{betasj}, $\Gamma^c(t)$ given by \eqref{gmmt}  satisfies \begin{equation}\label{ineq}
 \begin{split}
    \Gamma^c(t)\leq \frac{1}{a}\Big((a+\gamma\rho)d^{2}(nh)&e^{a(t-nh)}-\gamma\rho d^2(nh)\\&-\gamma a m(nh)\Big)e^{-\eta(t-nh)},
\end{split}
\end{equation}
where $a$ is given by \eqref{hhnm}, and $h$ is the sampling period chosen as in \eqref{hjgf}, along the solution of  \eqref{ctpe1}-\eqref{aqwr},\eqref{etcl},\eqref{dt},\eqref{obetbc3m},\eqref{petg}-\eqref{tildeGqwref} for all $t\in\big[nh,(n+1)h\big)$ and any $n\in\big[t^p_j/h,t^p_{j+1}/h\big)\cap\mathbb{N}$. 
\end{lem}
\noindent\textbf{Proof of Lemma \ref{imptnt}.} Taking the time derivative of \eqref{gmmt} in $t\in(nh,(n+1)h)$ and $n\in\big[t_j^p/h,t^p_{j+1}/h\big)\cap\mathbb{N}$, using Young's inequality, the relation  \eqref{ghm}, the dynamics of $m(t)$ given by \eqref{obetbc3m}, and \eqref{gmmt} to substitute for  $d^2(t)$, we show that 
\begin{equation*}
\begin{split}
    \dot{\Gamma}^c(t )\leq&(1+\rho_1+\gamma\rho)\Gamma^{c}(t)+\gamma(a+\gamma\rho)m(t)\\&-(\gamma\beta_1-\alpha_1)\Vert \hat{u}[t]\Vert^2-(\gamma\beta_2-\alpha_2)\hat{u}^2(1,t)\\&-(\gamma\beta_3-\alpha_3)\big(\theta_1\tilde{u}^2(0,t)+\theta_2\tilde{u}^2(1,t)\big).
\end{split}
\end{equation*}
Noting that both sides of this inequality are well-behaved in $t\in(nh,(n+1)h)$ and $n\in\big[t^p_j/h,t^p_{j+1}/h\big)\cap\mathbb{N}$, we assert, there exists a non-negative function $\iota(t )\in C^0\big((t^p_{j},t^p_{j+1});\mathbb{R}_{+}\big)$ such that 
\begin{equation}\label{zmkl}
\begin{split}
    \dot{\Gamma}^c(t )=&(1+\rho_1+\gamma\rho)\Gamma^{c}(t)+\gamma(a+\gamma\rho)m(t)\\&-(\gamma\beta_1-\alpha_1)\Vert \hat{u}[t]\Vert^2-(\gamma\beta_2-\alpha_2)\hat{u}^2(1,t)\\&-(\gamma\beta_3-\alpha_3)\big(\theta_1\tilde{u}^2(0,t)+\theta_2\tilde{u}^2(1,t)\big)-\iota(t),
\end{split}
\end{equation}
 for all $t\in(nh,(n+1)h)$ and $n\in\big[t^p_j/h,t^p_{j+1}/h\big)\cap\mathbb{N}$. 
Furthermore, using \eqref{gmmt} to substitute for  $d^2(t)$, we rewrite the dynamics of $m(t)$ as 
\begin{equation}\label{appop}
\begin{split}
    \dot{m}(t) =& -\rho \Gamma^{c}(t )-(\gamma\rho+\eta)m(t)+\beta_1\Vert \hat{u}[t]\Vert^2\\&+\beta_2 \hat{u}^2(1,t)+\beta_3\big(\theta_1\tilde{u}^2(0,t)+\theta_2\tilde{u}^2(1,t)\big),
\end{split}
\end{equation}
for $t\in(nh,(n+1)h)$ and $n\in\big[t^p_j/h,t^p_{j+1}/h\big)\cap\mathbb{N}$. Then, concatenating the equations \eqref{zmkl} and \eqref{appop}, we obtain the following ODE system 
\begin{equation}\label{zv}
    \dot{z}(t )=Az(t )+v(t),
\end{equation}where
\begin{equation*}
\begin{split}
&z(t )=
\begin{bmatrix}
\Gamma^{c}(t )\\m(t)
\end{bmatrix},
\text { }
    A = \begin{bmatrix}
1+\rho_1+\gamma\rho &  \gamma\big(a+\gamma\rho\big)\\
    -\rho & -(\gamma\rho+\eta)
    \end{bmatrix},\\&
    v(t )  =\begin{bmatrix}
    \Big(
    \begin{split}
        &-(\gamma\beta_1-\alpha_1)\Vert \hat{u}[t]\Vert^2-(\gamma\beta_2-\alpha_2)\hat{u}^2(1,t)\\&-(\gamma\beta_3-\alpha_3)\big(\theta_1\tilde{u}^2(0,t)+\theta_2\tilde{u}^2(1,t)\big)-\iota(t )
        \end{split}\Big)\\\begin{split}\Big(\beta_1\Vert \hat{u}[t]\Vert^2&+\beta_2\hat{u}^2(1,t)\\&+\beta_3\big(\theta_1\tilde{u}^2(0,t)+\theta_2\tilde{u}^2(1,t)\big)\Big)\end{split}
    \end{bmatrix}.
\end{split}
\end{equation*}
From \eqref{zv}, we obtain that
\begin{equation*}
    z(t )= e^{A(t-nh)}z(nh )+\int_{nh}^{t} e^{A(t-\xi)}v(\xi )d\xi,
\end{equation*}
for all $t\in[nh,(n+1)h)$ and $n\in\big[t^p_j/h,t^p_{j+1}/h\big)\cap\mathbb{N}$, using which we obtain
\begin{equation*}
    \Gamma^{c}(t )= Ce^{A(t-nh)}z(nh )+\int_{nh}^{t} Ce^{A(t-\xi)}v(\xi )d\xi,
\end{equation*}
where $C=\begin{bmatrix}
        1&0
    \end{bmatrix}$. The matrix $A$ has two distinct eigenvalues $-\eta$ and $1+\rho_1$. Therefore, using the matrix diagonalization of $A$ and after some simplifications, it can be shown that
\begin{equation*}
\begin{split}
    &Ce^{A(t-\xi)}v(\xi )=-\big((\gamma\beta_1-\alpha_1)g_1(t-\xi)-\beta_1 g_2(t-\xi)\big)\Vert \hat{u}[\xi]\Vert^2\\&\quad\qquad\qquad-\big((\gamma\beta_2-\alpha_2)g_1(t-\xi)-\beta_2 g_2(t-\xi)\big)\hat{u}^2(1,\xi)\\&\quad\qquad\qquad-\theta_1\big((\gamma\beta_3-\alpha_3)g_1(t-\xi)-\beta_3 g_2(t-\xi)\big)\tilde{u}^2(0,t)\\&
    \quad\qquad\qquad-\theta_2\big((\gamma\beta_3-\alpha_3)g_1(t-\xi)-\beta_3 g_2(t-\xi)\big)\tilde{u}^2(1,t)\\&\quad\qquad\qquad -g_1(t-\xi)\iota(\xi),
\end{split}
\end{equation*}
where
\begin{equation*}
    g_1(t)=\frac{1}{a}\big(-\gamma\rho +(a+\gamma\rho)e^{at}\big)e^{-\eta t},
\end{equation*}
and
\begin{equation*}
    g_2(t)=\frac{\gamma (a+\gamma\rho)}{a}\big(-1+ e^{at}\big)e^{-\eta t}.
\end{equation*}
We can easily observe that $g_1(t)>0$ for all $t\geq 0$. Furthermore, noting that $\gamma\beta_i/\alpha_i=1/(1-\sigma), i=1,2,3$ from \eqref{betasj}, and recalling \eqref{mdt}, we show that
\begin{equation*}
\begin{split}
    &\big(\gamma\beta_i-\alpha_i\big)g_1(t-\xi)-\beta_i g_2(t-\xi)\\&=\frac{\alpha_i (a+\gamma\rho) }{a}\Big(1+\frac{\sigma a}{(1-\sigma)(a+\gamma\rho)}-e^{a(t-\xi)}\Big)e^{-\eta (t-\xi)}\\&=\frac{\alpha_i (a+\gamma\rho) }{a}\Big(e^{a\tau}-e^{a(t-\xi)}\Big)e^{-\eta (t-\xi)},
\end{split}
\end{equation*}
for $i=1,2,3$. As $nh\leq\xi\leq t<(n+1)h$, and $h\leq \tau$, we have that $(\gamma\beta_i-\alpha_i)g_1(t-\xi)-\beta_i g_2(t-\xi)>0$ for $i=1,2,3$. Thus, we argue that $Ce^{A(t-\xi)}v(\xi )\leq 0$ for all $t,\xi$ such that $nh\leq\xi\leq t<(n+1)h$, and $n\in\big[t^p_j/h,t^p_{j+1}/h\big)\cap\mathbb{N}$. Considering this, it can be derived for $t\in[nh,(n+1)h)$ that 
\begin{equation*}
\begin{split}
    &\Gamma^{c}(t )\leq Ce^{A(t-nh)}z(nh )\\&\leq g_1(t-nh)\Gamma^{c}(nh )+g_2(t-nh)m(nh)\\&\leq 
\frac{1}{a}\Big(-\gamma(a+\gamma\rho)m(nh)-\gamma\rho\Gamma^{c}(nh )\\&\quad\qquad+(a+\gamma\rho)\big(\Gamma^{c}(nh )+\gamma m(nh)\big)e^{a(t-nh)}\Big)e^{-\eta(t-nh)}.
\end{split}
\end{equation*}
By substituting for $\Gamma^{c}(nh )$ using \eqref{gmmt}, we obtain the inequality \eqref{ineq} that is valid for $t\in[nh,(n+1)h)$. This completes the proof of Lemma \ref{imptnt} \hfill$\qed$

Now let us continue with the proof of Theorem \ref{hhgbnml}. Assume that an event has triggered at $t=t^p_j$ and $m(t^p_j)>0$. Then, let us analyze the behavior of $\Gamma^c(t)$ and $m(t)$ in $t\in[t^p_j,t^p_{j+1})$ along the solution of \eqref{ctpe1}-\eqref{aqwr},\eqref{etcl},\eqref{dt},\eqref{obetbc3m},\eqref{petg}-\eqref{tildeGqwref}. After the event at $t=t_j^p$, the control law is updated. Thus, we have from \eqref{gmmt} that $\Gamma^c(t^p_j)=-\gamma m(t_j^p)<0$. Consequently, $\Gamma^c(t)$ will at least remain non-positive until $t=t^p_j+\tau$ where $\tau$ is the MDT given by \eqref{mdt} (see R1 of Theorem \ref{hhgb}). Thus, $\Gamma^c(t)$ will definitely remain non-positive in $t\in[t^p_j,t^p_j+h)$ as $h\leq\tau$. However, at each $t=nh,n>0,$ the periodic event-trigger given by \eqref{petg}-\eqref{tildeGqwref} is evaluated, and only if $\Gamma^p(nh)>0$ that an event is triggered, and the control input is updated. If $\Gamma^p(nh)\leq 0$, then an update would not be required as $\Gamma^c(t)$ will be non-positive due to the relation \eqref{ineq} (Note that the right hand side of \eqref{ineq} is definitely non-positive when $\Gamma^p(nh)\leq 0$). Thus, $\Gamma^c(t)$ will in fact remain non-positive at least until $t=t_{j+1}^p$ where $\Gamma^p(t_{j+1}^{p-})> 0$. As $\Gamma^c(t)\leq 0$ for $t\in[t_{j}^p,t_{j+1}^p)$, we write from \eqref{gmmt} that $d^2(t)\leq \gamma m(t)$ for $t\in[t_j^p,t_{j+1}^p)$. Then, considering the dynamics of $m(t)$ given by \eqref{obetbc3m}, we get  $\dot{m}(t)\geq -(\eta+\gamma\rho)m(t)$ for $t\in(t^p_j,t^p_{j+1})$, which leads to $m(t)\geq e^{-(\eta+\gamma\rho)(t-t_j^p)}m(t_j^p)>0$ for $t\in[t_j^p,t_{j+1}^p)$. The time continuity of $m(t)$ leads to $m(t_{j+1}^{p-})=m(t_{j+1}^p)>0$. Therefore, after the control input has been updated at $t=t_{j+1}^p$, we obtain the equality $\Gamma^c(t_{j+1}^p)=-\gamma m(t^p_{j+1})<0$. In a similar way, we can analyze the behavior of $\Gamma^c(t)$ and $m(t)$ in all $t\in[t^p_j,t^p_{j+1})$ for any $j\in\mathbb{N}$ starting from the first event at $t^p_0=0$ where  $m(0)>0$ to prove that $\Gamma^c(t)\leq 0$ for all $t\in[t^p_j,t^p_{j+1}),j\in\mathbb{N}$ and $m(t)>0$ for all $t>0$. Thus, the global $L^2$-exponential convergence of the closed-loop system to zero satisfying the estimate \eqref{rrt} follows from R4 of Theorem \ref{hhgb}. This completes the proof of Theorem \ref{hhgbnml}. \hfill $\qed$

\medskip

\noindent \textbf{B. Proof of Theorem \ref{hhgbsd}}

The well-posedness of the closed-loop system \eqref{ctpe1}-\eqref{ctoe3kl} with $\theta_1=0,\theta_2=1$ under the STC is discussed in Remark \ref{ici}. To streamline the proof of Theorem \ref{hhgbsd}, we first present Lemma \ref{zbmmk}.  

\begin{lem}\label{zbmmk} Consider the STC approach \eqref{xxvb}-\eqref{dsfg} under Assumptions \ref{ass1} and \ref{ddcf}, which generates an increasing set of event times $\{t^s_{j}\}_{j\in\mathbb{N}}$ with $t^s_{j}=0$. Then, for the error $d(t)$ given by  \eqref{dt} and $m(t)$ governed by \eqref{obetbc3m}, the followings hold:
\begin{equation}\label{bbv1aq}
    d^2(t)\leq H(t^s_j)e^{2\varrho(t-t^s_j)},
\end{equation}
and
\begin{equation}\label{bbv2aq}
   m(t)\geq m(t^s_j)e^{-\eta(t-t^s_j)}-\frac{\rho H(t^s_j)}{2\varrho+\eta}e^{-\eta (t-t^s_j)}\big(e^{(2\varrho+\eta)(t-t^s_j)}-1\big),
\end{equation}
for all $t\in[t^s_{j},t^s_{j+1}),j\in\mathbb{N}$ along the solution of \eqref{ctpe1}-\eqref{aqwr}, where  $H(t)$ and $\varrho$  are given by \eqref{xxcvb} and \eqref{qqqwwe}, respectively. 
\end{lem}
\noindent\textbf{Proof of Lemma \ref{zbmmk}.} Consider the positive definite function
\begin{equation}\label{fffg}
    V = \frac{1}{2}\int_{0}^{1}\hat{u}^2(x,t)dx.
\end{equation}
Taking its time derivative  along the solution \eqref{ctoe1kl}-\eqref{ctoe3kl} and using Young's inequality and Cauchy-Schwarz inequality, we show that
\begin{equation*}\label{ddsa}
    \begin{split}
        \dot{V}\leq &-\varepsilon q\hat{u}^2(1,t)-\varepsilon\Vert\hat{u}_x[t]\Vert^2+\lambda\Vert\hat{u}[t]\Vert^2+\frac{\varepsilon h_1}{2}\hat{u}^2(1,t)\\&+\frac{\varepsilon}{2h_1}(U_j^s)^2+\frac{\varepsilon p_{10}}{2h_2}\hat{u}^{2}(1,t)+\frac{\varepsilon p_{10}h_2}{2}\tilde{u}^{2}(1,t)\\&+\frac{1}{2h_3}\tilde{u}^2(1,t)+\frac{\Vert p_1\Vert^2h_3}{2}\Vert\hat{u}[t]\Vert^2,
    \end{split}
\end{equation*}
for $t\in(t_{j}^s,t^s_{j+1}),j\in\mathbb{N}$ and some $h_1,h_2,h_3>0$. Let us select $h_1=\frac{\lambda}{2\varepsilon},h_2=\frac{2\varepsilon p_{10}}{\lambda},h_3=1$. Then, one can show
\begin{equation}\label{sssf}
\begin{split}
    \dot{V}\leq& -\varepsilon\Big( q-\frac{\lambda}{2\varepsilon}\Big)\hat{u}^2(1,t)-\varepsilon\Vert\hat{u}_x[t]\Vert^2+\varrho\Vert\hat{u}[t]\Vert^2\\&+\frac{\varepsilon^2}{\lambda}(U^s_j)^2
+\Big(\frac{\varepsilon^2p_{10}^2}{\lambda}+\frac{1}{2}\Big)\tilde{u}^2(1,t),
\end{split}
\end{equation}
for $t\in(t^s_{j},t^s_{j+1}),j\in\mathbb{N}$ where $\varrho>0$ is given by \eqref{qqqwwe}. Using Cauchy-Schwarz inequality and \eqref{etcl} under the consideration of \eqref{fffg}, we obtain that $(U^s_j)^2\leq 2\Vert k\Vert^2 V(t^s_j).$ Thus, recalling Assumption \ref{ass1} from which it follows that $q>\lambda/2\varepsilon$ for $\theta_1=0$, we write \eqref{sssf} as
\begin{equation}\label{dddsa}
\begin{split}
    &\dot{V}\leq 2\varrho V(t)+\frac{2\varepsilon^2\Vert k\Vert^2}{\lambda}V(t^s_j)+\Big(\frac{\varepsilon^2p_{10}^2}{\lambda}+\frac{1}{2}\Big)\tilde{u}^{2}(1,t),
\end{split}
\end{equation}
for $t\in(t^s_{j},t^s_{j+1}),j\in\mathbb{N}$. Now let us find a known upper-bound for $\vert\tilde{u}(1,t)\vert$. From Lemma 1 of \cite{rathnayake2022sampled}, we obtain that $\Vert\tilde{w}_{x}[t]\Vert\leq \big(\Vert\tilde{w}_{x}[0]\Vert+M_{1}\Vert\tilde{w}[0]\Vert\big)e^{-\sigma^*t}$ and $\Vert\tilde{w}[t]\Vert\leq \Vert\tilde{w}[0]\Vert e^{-\sigma^*t}$ for all $t\geq 0$ where $\tilde{w}$ is the observer error target state obtained via \eqref{qmpo} with $\theta_1=0,\theta_2=1$ (see \cite{rathnayake2022sampled} for details), $M_1$ is given by \eqref{phmk1}, and $\sigma^*$ is given by \eqref{qqqwwe}. Thus, using \eqref{qppu}, Agmon's and Young's inequalities, we obtain that
\begin{equation}\label{zznnm}
\begin{split}
    &\vert \tilde{u}(1,t)\vert=\vert \tilde{w}(1,t)\vert\leq (\sqrt{2})^{-1}\Vert\tilde{w}[t]\Vert+(\sqrt{2})^{-1}\Vert\tilde{w}_x[t]\Vert\\&\leq (\sqrt{2})^{-1}(M_1+1)\Vert\tilde{w}[0]\Vert e^{-\sigma^*t}+(\sqrt{2})^{-1}\Vert\tilde{w}_x[0]\Vert e^{-\sigma^*t},
\end{split}
\end{equation}
for all $t\geq 0$. But, using \eqref{qmpo} with $\theta_1=0,\theta_2=1$, and Cauchy-Schwarz inequality, we show that $\Vert\tilde{w}[0]\Vert\leq \Omega_1\Vert \tilde{u}[0]\Vert$ and $\Vert\tilde{w}_x[0]\Vert\leq \Vert\tilde{u}_x[0]\Vert+\Omega_2\Vert\tilde{u}[0]\Vert$ where $\Omega_1$ and $\Omega_2$ are given by \eqref{zzzs1} and \eqref{dsfg}, respectively. Again using Cauchy-Schwarz inequality, we obtain from \eqref{aqwr} that $\Vert\tilde{u}[0]\Vert\leq \Vert u[0]\Vert+\Vert\hat{u}[0]\Vert$ and $ \Vert\tilde{u}_x[0]\Vert\leq \Vert u_x[0]\Vert+\Vert\hat{u}_x[0]\Vert.$ Thus, from \eqref{zznnm}, we obtain the following bound for $\tilde{u}^2(1,t)$
\begin{equation}\label{aaaq1}
    \vert\tilde{u}(1,t)\vert\leq \Psi_0 e^{-\sigma^*t},
\end{equation}
for all $t\geq 0$ where $\Psi_0=(\sqrt{2})^{-1}\big((M_1+1)\Omega_1+\Omega_2\big)\Psi_1+(\sqrt{2})^{-1}\Psi_2+(\sqrt{2})^{-1}\big((M_1+1)\Omega_1+\Omega_2\big)\Vert\hat{u}[0]\Vert+(\sqrt{2})^{-1}\Vert\hat{u}_x[0]\Vert.$ Then, considering \eqref{aaaq1} and noting that $e^{-2\sigma^*t}<e^{-2\sigma^*t_j^s}$ for all $t>t_j^s$, we obtain from \eqref{dddsa} that 
\begin{equation}\label{sssffg}
    \dot{V}(t)\leq 2\varrho V(t)+\frac{2\varepsilon^2\Vert k\Vert^2}{\lambda}V(t^s_j)+(\Psi^*_0)^2e^{-2\sigma^*t_j^s},
\end{equation}
for $t\in(t^s_{j},t^s_{j+1}),j\in\mathbb{N}$ where $\Psi_0^*$ is given by \eqref{ytre}. Therefore, from \eqref{sssffg}, we show that 
\begin{equation*}
\begin{split}
    V(t)\leq &e^{2\varrho(t-t^s_j)}V(t_j^s)\\&+\frac{\frac{2\varepsilon^2\Vert k\Vert^2}{\lambda}V(t^s_j)+(\Psi^*_0)^2e^{-2\sigma^*t_j^s}}{2\varrho}\big(e^{2\varrho(t-t^s_j)}-1\big),
\end{split}
\end{equation*}
for $t\in[t^s_j,t^s_{j+1}),j\in\mathbb{N}$ from which we obtain that
\begin{equation}\label{dddaxaq}
\begin{split}
    &\Vert \hat{u}[t]\Vert^2\\&\leq \Big(\Vert\hat{u}[t^s_j]\Vert^2+\frac{\varepsilon^2\Vert k\Vert^2}{\lambda\varrho}\Vert\hat{u}[t^s_j]\Vert^2+\frac{(\Psi_0^*)^2e^{-2\sigma^*t_j^s}}{\varrho}\Big)e^{2\varrho(t-t^s_j)},
\end{split}
\end{equation}
considering \eqref{fffg}. Using Cauchy-Schwarz inequality and Young's inequality on \eqref{dt}, we show that $ d^2(t)\leq 2\Vert k\Vert^2\Vert\hat{u}[t^s_j]\Vert^2+2\Vert k\Vert^2\Vert\hat{u}[t]\Vert^2.$ Then, using \eqref{dddaxaq}, we obtain \eqref{bbv1aq}. Considering the dynamics of $m(t)$ given by \eqref{obetbc3m} and the relation \eqref{bbv1aq}, we show $\dot{m}(t)\geq -\eta m(t)-\rho H(t^s_j)e^{2\varrho(t-t^s_j)}$ for $t\in(t^s_{j},t^s_{j+1}),j\in\mathbb{N}$ from which we obtain \eqref{bbv2aq}. This completes the proof of Lemma \ref{zbmmk}. \hfill $\qed$

Now let us continue with the proof of Theorem \ref{hhgbsd}. Consider the triggering function $\Gamma^c(t)$ given by \eqref{gmmt} along the solution of \eqref{ctpe1}-\eqref{aqwr},\eqref{etcl},\eqref{dt},\eqref{obetbc3m},\eqref{xxvb}-\eqref{dsfg}. Further, let us assume that an event has occurred at $t=t_j^s$ and $m(t_j^s)>0$. Then, as the control input is updated, it follows from \eqref{gmmt} that $\Gamma^c(t_j^s)=-\gamma m(t_j^s)<0$. Moreover, $\Gamma^c(t)$ will remain non-positive at least until $t=t_j^s+\tau$, where $\tau$ is the MDT given by \eqref{mdt} (see R1 of Theorem \ref{hhgb}). We have from \eqref{bbv1aq} that $d^2(t)\leq H(t_j^s)e^{2\varrho(t-t_j^s)}$ and from \eqref{bbv2aq} that 
\begin{equation*}\label{bbv211q}
\begin{split}
    \gamma m(t)\geq &\gamma m(t_j^s)e^{-\eta(t-t_j^s)}\\&-\frac{\gamma\rho H(t_j^s)}{2\varrho+\eta}e^{-\eta (t-t_j^s)}\big(e^{(2\varrho+\eta)(t-t_j^s)}-1\big),
\end{split}
\end{equation*}
for $t\in[t_j^s,t_{j+1}^s)$. Note that the RHS of \eqref{bbv1aq}
is an increasing function of $t$ whereas the RHS of \eqref{bbv2aq} is a decreasing function of $t$. Then, if there is a positive $t^\dagger>t_j^s$ that satisfies   
\begin{equation}\label{alfg}
\begin{split}
H(t_j^s)e^{2\varrho(t^\dagger-t_j^s)}&= \gamma m(t_j^s)e^{-\eta(t^\dagger-t_j^s)}\\&-\frac{\gamma\rho H(t_j^s)}{2\varrho+\eta}e^{-\eta (t^\dagger-t_j^s)}\Big(e^{(2\varrho+\eta)(t^\dagger-t_j^s)}-1\Big),
\end{split}
\end{equation}
we are certain that $d^2(t)\leq \gamma m(t)$, \textit{i.e.,} $\Gamma^c(t)\leq 0$ for $t\in[t_j^s,t^\dagger)$ (note that the LHS of \eqref{alfg} is an upper-bound for $d^2(t^\dagger)$, and the RHS of \eqref{alfg} is a lower-bound for $\gamma m(t^\dagger)$). Solving \eqref{alfg} for $t^\dagger$, we obtain that 
\begin{equation*}
t^\dagger = t_j^s+\frac{1}{2\varrho+\eta}\ln\left(\frac{\gamma m(t_j^s)+\frac{\gamma\rho H(t_j^s)}{2\varrho+\eta}}{H(t_j^s)+\frac{\gamma\rho H(t_j^s)}{2\varrho+\eta}}\right).
\end{equation*}
If $t^\dagger>t_j^s+\tau$, the next event can be chosen as $t_{j+1}^s=t^\dagger$. If $t^\dagger\leq t_j^s+\tau$, the next event can be chosen as $t_{j+1}^s=t_j^s+\tau$. In this way, since the next event time is given by \eqref{xxvb}-\eqref{dsfg}, it is ensured that $\Gamma^c(t)\leq 0$ for $t\in [t^s_j,t^s_{j+1})$ while preventing the occurrence of Zeno phenomenon. Then, employing the same line of reasoning as in the proof of Theorem \ref{hhgbnml}, we can show that $\Gamma^c(t)\leq 0$ for all $t\in[t^s_j,t^s_{j+1}),j\in\mathbb{N}$ and $m(t)>0$ for all $t>0$. Thus, the $L^2$-exponential convergence to zero satisfying the estimate \eqref{rrt} follows from R4 of Theorem \ref{hhgb}. This completes the proof of Theorem \ref{hhgbsd}.  \hfill $\qed$
   
\bibliographystyle{IEEEtranS}
\bibliography{main}

\end{document}